\documentclass{article}

\usepackage{arxiv}

\usepackage[utf8]{inputenc} % allow utf-8 input
\usepackage[T1]{fontenc}    % use 8-bit T1 fonts
\usepackage{hyperref}       % hyperlinks
\usepackage{booktabs}       % professional-quality tables
\usepackage[sort]{natbib}
\usepackage{doi}
\usepackage{subfig}
\usepackage[inline]{enumitem}

\usepackage{pgf,tikz}
\usetikzlibrary{matrix,arrows,automata,positioning,arrows.meta,calc}

\usepackage{amsmath,amssymb,mathtools}

\usepackage{amsthm}

\newtheorem{lemma}{Lemma}
\newtheorem{theorem}{Theorem}
\newtheorem{corollary}{Corollary}
\newtheorem{conjecture}{Conjecture}
\newtheorem{helpful}{Helpful Fact \arabic{helpful}}

\definecolor{pureBlack}{rgb}{0,0,0}

\usepackage[mathlines]{lineno}
\allowdisplaybreaks[2]

\def\VFirst {{"1","2","3","4","5","6","7","8","9","10","11","12","13","14","15","16"}}
\def\VSecond{{"12","14","11","6","15","13","5","1","10","8","2","16","3","4","7","9"}}
\def\VThird {{"13","8","12","6","7","15","4","11","16","1","9","5","3","10","2","14"}}

\title{Valid path-based graph vertex numbering}

\date{}

\author{ \href{https://orcid.org/0000-0000-0000-0000}{\includegraphics[scale=0.06]{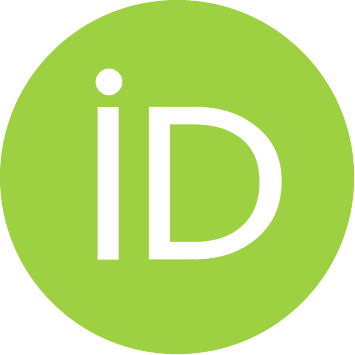}\hspace{1mm}Les Foulds}
 \quad
	\href{https://orcid.org/0000-0000-0000-0000}{\includegraphics[scale=0.06]{orcid.pdf}\hspace{1mm}Humberto J. Longo} \\
 Instituto de Inform\'{a}tica, Universidade Federal de Goi\'{a}s,\\
 Alameda Palmeiras, Quadra D, Campus Samambaia,\\
 CEP 74690-900, Goi\^{a}nia, Goi\'{a}s, Brazil\\
	\texttt{lesfoulds,longo@ufg.br} \\
}

\renewcommand{\shorttitle}{Valid path-based graph vertex numbering}

\hypersetup{
pdftitle={Valid path-based graph vertex numbering},
pdfsubject={},
pdfauthor={Les Foulds, Humberto J. Longo},
pdfkeywords={Graph, numbering, 2-path, validity, $\mathcal{NP}$-hard},
}

\begin{document}
% \linenumbers
\maketitle

\begin{abstract}
	A labelling of a graph is an assignment of labels to its vertex or edge sets (or both), subject to certain conditions, a well established concept. A labelling of a graph $G$ of order $n$ is termed a \emph{numbering} when the set of integers $\{1,\dots,n\}$ is used to label the vertices of $G$ distinctly. A 2-path (a path with three vertices) in a vertex-numbered graph is said to be \emph{valid} if the number of its middle vertex is smaller than the numbers of its endpoints. The problem of finding a vertex numbering of a given graph that optimises the number of induced valid 2-paths is studied, which is conjectured to be in the $\mathcal{NP}$-hard class. The reported results for several graph classes show that apparently there are not one or more numbering patterns applicable to different classes of graphs, which requires the development of a specific numbering for each graph class under study.
\end{abstract}

% keywords can be removed

%--------------------------------------------------- S E C T I O N -
\section{Introduction}
\label{introduction}
A labelling of a graph is an assignment of labels to its vertex or edge sets (or both), subject to certain conditions, a concept first introduced by \citet{Rosa-1967}. Since that seminal work, thousands of papers on graph labelling have appeared, as dynamically surveyed by  \citet{Gallian-2018}. When used without qualification, the term \textit{labelled graph} generally refers to a vertex-labelled graph. Several well-known optimisation problems on graphs can be formulated in terms of labelling as follows: given a graph $G$, find a labelling of $G$ such that a given linear function of the chosen labelling is optimised. Applications fields include, among many others, network addressing \citep{Judd-Beer-1995}, circuit layout code design \citep{Shieh-Yu-Yen-Yu-Chang-Xu-Wann-Taiwan-2011} and the bandwidth, minimum-sum and cutwidth problems \citep{Chung-1981,Juvan-Mohar-1992}. \emph{Network addressing} involves messages that are sent from source nodes to destination nodes that define paths over which the messages should travel to reach their destination nodes. \emph{A circuit layout code} represents the layout of an integrated circuit device. The \emph{bandwidth problem} involves finding the labelling that minimises the maximum stretch over all the edges. The \emph{minimum sum problem} involves finding a labelling that minimises the total length sum over the edges.  The \emph{cutwidth problem} involves finding a labelling that minimises the maximum overlap over the edges.  

A labelling of a graph $G$ of order $n$ is termed a \emph{numbering} (\emph{ordering}) when the set of integers $\{1,\dots,n\}$ is used to label the vertices of $G$ distinctly. In this case, the vertices can be arranged along a line or a linear array. Numbered graphs were popularised by \citet{Bloom-Golomb-1977} who observed that they are useful as mathematical models for a broad range of applications, such as the analysis of various coding theory problems, including the design of good radar-type codes, sync-set and convolutional codes with optimal autocorrelation properties. They also mentioned applications to problems in additive number theory, determination of ambiguities in X-ray crystallographic analysis, design of communication network addressing systems and identification of optimal circuit layouts.

In this paper the concept of \emph{path validity} is introduced. A path of a numbered graph is said to be \emph{valid} if it is a 2-path (a path with three vertices) and the number of its middle vertex is strictly smaller than the numbers of its endpoints. Note that the concept of validity is defined only for 2-paths, not for longer paths. The \emph{valid path problem} introduced here involves finding a numbering of a given graph $G$ that optimises the quantity of valid paths in $G$ with respect to the numbering. Results for various classes of graphs are reported.

This problem arises, for example, in the enumeration of chordless cycles in graphs, which has application in the prediction of nuclear magnetic resonance, chemical shift values and ecological networks \citep{Dias-Castonguay-Longo-Jradi-2014}. Furthermore, the valid path problem has further applications, concerning numbered graphs that are models of various discrete systems in which only partial comparisons of the vertex numbers can be made. One such case occurs when, due to the correlation structure of the system, the only vertex number comparison possible involves identifying the 2-paths that have a middle vertex with the lowest number (have numbers of the form 2,1,3).

Here, the aim is to make an assignment that optimises the number of such 2-paths. Such models arise in the analysis of financial data -- where the vertices represent monetary entities such as stocks or shares and the edges represent correlations between them \citep{Massara-DiMatteo-Tomaso-2016}, facility (plant) layout -- where the vertices represent facility objects and the edges represent the adjacencies between them in a plan of the facility \citep{Seppanen-Moore-1970}, integrated circuit design -- where the vertices represent the electrical elements of a circuit and the edges represent the physical connections between them \citep{Lengauer-2012}, systems biology -- where the vertices represent proteins and the edges represent protein interactions in a metabolic network \citep{Song-Aste-DiMatteo-2008} and  social systems -- where the vertices represent social agents (e.g. individuals, groups or companies)  and the edges represent social interactions \citep{Easley-Kleinberg-2010}.

%--------------------------------------------------- S E C T I O N -
\section{A formal definition of the path validity problem}
\label{sec:details}

Let $G=(V(G),E(G))$ be a finite, undirected, simple graph with vertex set $V(G)$ and edge set $E(G)$, where $n=|V(G)|$ and $m=|E(G)|$, and $Adj(u) = \{v \in V(G)\mid \{u, v\} \in E(G)\}$. A 2-path $\langle x,u,y\rangle$ exists in $G$ for every vertex $u\in V(G)$ such that $\{x,u\}, \{u,y\} \in E(G)$. The square of the adjacency matrix of $G$ counts the number of such paths \cite{Bondy-Murty-2008,Weisstein-2017}.

A numbering $\pi$ of $G$ is a bijective mapping from $V(G)$ to $\{1, 2,\dots ,n\}$. If $G$ has a such numbering, a 2-path $\langle x, u, y\rangle$ in $G$ is said to be valid if the number associated with its middle vertex $u$ is smaller than the numbers of its endpoints $x$ and $y$, i.e $\pi(u) <\min\{\pi(x),\pi(y)\}$. The squared adjacency matrix can be used to check if there is the possibility of a valid path existing between any given pair of vertices.

The validity $\phi_\pi(G)$ of a numbering $\pi$ of $G$ is defined as the quantity of valid 2-paths of $G$ induced by $\pi$, i.e.,
\begin{equation}\label{eq:validity}
\!\phi_\pi(G) = |\{\langle x, u, y\rangle\mid u \in V(G),\> x, y \in Adj(u)\text{ and } \pi(u) < \min(\pi(x), \pi(y))\}|.\!
\end{equation}

The \emph{minimum validity} $\phi_{\min}(G)$ of $G$ is the minimum of $\phi_\pi(G)$ over all numberings $\pi$, i.e.,
\begin{equation}\label{eq:minval_g}
\phi_{\min}(G) = \min\{\phi_\pi(G)\mid\pi\text{ is a numbering of }G\}.
\end{equation}
Analogously, the \emph{maximum validity} $\phi_{\max}(G)$ of $G$ is the maximum of $\phi_\pi(G)$, i.e.,
\begin{equation}\label{eq:maxval_g}
\phi_{\max}(G) = \max\{\phi_{\pi}(G)\mid \pi\text{ is a numbering of }G\}.
\end{equation}

For general graphs of nontrivial order, $\phi_{\min}(G)$ and $\phi_{\max}(G)$ can vary significantly, depending on the graph structure and how the numbering is carried out. If no constraints are imposed then, of course, $n!$ different numberings are possible, inducing sets of valid paths that may vary quite markedly in cardinality and element composition.
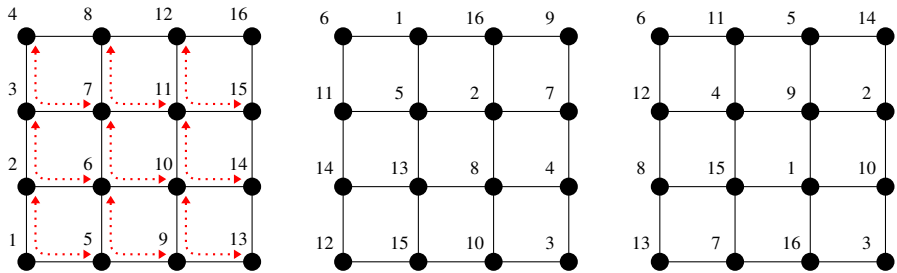
\begin{figure}[htpb!]
 \centering
 \begin{tikzpicture}[
       pil/.style={
           {Triangle[scale=0.6]}-{Triangle[scale=.6]},
           thick,
           dotted,
           rounded corners,
           red
           },
   scale=2,transform shape]
  \matrix[nodes={ultra thick},row sep=0.3cm,column sep=0.6cm]
  {
   \foreach \x in {0,1,2,3} {
    \foreach \y in {0,1,2,3} {
     \node[font={\scriptsize},align=flush  right] (a) at (\x-0.18,\y+0.3) {\pgfmathparse{\VFirst[\x*4+\y])}\pgfmathresult};
     \fill[color=black] (\x,\y) circle (0.12);
    }
   }
   \foreach \x [evaluate=\x as \i using \x+1] in {0,1,2} {
    \foreach \y [evaluate=\y as \j using \y-1] in {3,2,1} {
     \draw [pil] ($(\x,\y)+ (.12,-.12)$) -- ($(\x,\j)+ (.12,.10)$) -- ($(\i,\j)+ (-.14,.10)$);
    }
   }
   \draw (0,0) grid (3,3); &

   \foreach \x in {3,2,1,0} {
    \foreach \y in {3,2,1,0} {
     \node[font={\scriptsize}] (a) at (\x-0.25,\y+0.25) {\pgfmathparse{\VSecond[\x*4+\y])}\pgfmathresult};
     \fill[color=black] (\x,\y) circle (0.12);   
    }
   }
   \draw (0,0) grid (3,3); &

   \foreach \x in {3,2,1,0} {
    \foreach \y in {0,1,2,3} {
     \node[font={\scriptsize}] (a) at (\x-0.25,\y+0.25) {\pgfmathparse{\VThird[\x*4+\y])}\pgfmathresult};
     \fill[color=black] (\x,\y) circle (0.12);   
    }
   }
   \draw (0,0) grid (3,3); \\
  };
 \end{tikzpicture}
 \caption{Three numberings of the $4\times 4$ grid graph.}
 \label{tikz:three_possible_labelling}
\end{figure}

Figure \ref{tikz:three_possible_labelling} illustrates this phenomenon for three labellings of the $4\times 4$ grid graph (a $p\times q$ grid graph is the graph Cartesian product $P_p\times P_q$ of the path graphs on $p$ and $q$ vertices). The leftmost numbering induces only the nine highlighted valid 2-paths (the minimum possible), the second has 15 of these structures and the rightmost has 26 valid 2-paths (the maximum possible). These three illustrative numberings were obtained with the aid of a brute force algorithm.

As far as we know, finding numberings that achieve $\phi_{\min}(G)$ and $\phi_{\max}(G)$ for an arbitrary graph $G$ are open problems, which are denoted by MIN-VP and MAX-VP, respectively, and are defined as:
\paragraph{\bf\underline{Instance}:} A finite, undirected, simple, unweighted graph $G$ and an integer $k$.
\paragraph{\bf\underline{MIN-VP}:}
\begin{description}[topsep=-2pt, itemsep=-2pt]
 \item [Question:] Is there a numbering $\pi$ of $G$ that induces at most $k$ valid paths in $G$, i.e., $\phi_\pi(G) \leqslant k$?
\end{description}
\paragraph{\bf \underline{MAX-VP}:}
\begin{description}[topsep=-2pt, itemsep=-2pt]
 \item [Question:] Is there a numbering $\pi$ of $G$ that induces at least $k$ valid paths in $G$, i.e., $\phi_\pi(G) \geqslant k$?
\end{description}

%--------------------------------------------------- S E C T I O N -
\section{Some results for particular classes of graphs}
\label{sec:some-results}
Some results are now stated for various classes of graphs. Because $\phi_\pi(G)$ is additive over disconnected components of $G$, we assume henceforth that $G$ is connected.

% __________________ P A T H S
\begin{lemma}\label{lbl-paths}
If $P_n$ is a path of finite order $n\geqslant 3$ (a linear graph) then: ($i$) $\phi_{\min}(P_n) = 0$ and ($ii$) $\phi_{\max}(P_n) = \lceil\frac{n}{2}\rceil-1$.
% \begin{align}
%   (i)\ \phi_{\min}(P_n) & = 0\label{eq:lemma1-1}\shortintertext{and}
%   (ii)\ \phi_{\max}(P_n) & = \left\lceil\frac{n}{2}\right\rceil-1.\label{eq:lemma1-2}
% \end{align}
\end{lemma}
\begin{proof} Let $u$ be a vertex at one end of the path $P_n$. 
\begin{enumerate}[leftmargin=*,label=(\emph{\roman*})]
 \item\label{lbl-path-i} Suppose that $u$ is numbered as 1. By sequentially numbering as $2, 3, \dots n$ the remaining vertices, from the vertex adjacent to $u$ to the opposite end of the path, it follows easily that no valid path exists. Thus, $\phi_{\min}(P_n)=0$. Furthermore, given this numbering, if two of its consecutive numbers, say $k$ and $k+1$, for $k=1,2,\dots n-2$, are interchanged, a valid path is produced and the new numbering is no longer minimal. Subsequent distinct pairwise number interchanges preserve the presence of a valid path unless the original numbering or its reverse is generated.
 \item\label{lbl-path-ii} Suppose that each vertex distant $2k+1$ from $u$, $k=0,\dots,\lfloor\frac{n}{2}\rfloor-1$, is numbered as $k+1$ and the remaining vertices are arbitrarily and distinctly numbered as $\ell\in\mathcal{L}=\{\lfloor\frac{n}{2}\rfloor+1,\dots,n\}$. If $n$ is even, each vertex numbered as $k+1$, for $k=0,\dots,\lfloor\frac{n}{2}\rfloor-2$, induces exactly one valid path, yielding $\lfloor\frac{n}{2}\rfloor-1=\lceil\frac{n}{2}\rceil-1$ valid paths in total. If $n$ is odd, each vertex numbered as $k+1$, for  $k=0,\dots,\lfloor\frac{n}{2}\rfloor-1$, induces exactly one valid path, yielding $\lfloor\frac{n}{2}\rfloor=\lceil\frac{n}{2}\rceil-1$ valid paths in total.

Given this numbering $\phi$ of $P_n$, suppose that the numbers $r_i$ and $r_j$, of any two vertices $i$ and $j$ of $P_n$, are exchanged giving a new numbering $\phi'$. If $r_i, r_j \in \mathcal{L}$ or $r_i, r_j \notin \mathcal{L}$, no new valid path is created by $\phi'$. If $r_i \notin \mathcal{L}$ and $r_j \in \mathcal{L}$, then $r_i$ induces the valid path $p_i=\langle i-1,i,i+1\rangle$, which is numbered as $\langle r_{i-1},r_i,r_{i+1}\rangle$. Since $r_j > r_i$, the new numbering $\langle r_{i-1},r_j,r_{i+1}\rangle$ of the path $p_i$ may not be valid. On the other hand, the vertex $j$ is the extremity of at most two valid paths ($\langle j-2,j-1,j\rangle$ and $\langle j,j+1,j+2\rangle$), numbered in $\phi$ as $\langle r_{j-2},r_{j-1},r_j\rangle$ and $\langle r_j,r_{j+1},r_{j+2}\rangle$, respectively. Since $r_i < r_j$, at most one of the numberings $\langle r_{j-2},r_{j-1},r_i\rangle$ and $\langle r_i,r_{j+1},r_{j+2}\rangle$ will be valid. Thus, no order relation between the numbers of $i$ and $j$ and the numbers of the other vertices in the listed 2-paths will increase the quantity of valid paths in $\phi'$. The case $r_i \in \mathcal{L}$ and $r_j \notin \mathcal{L}$ is similar.

Therefore, $\phi_{\max}=\lceil\frac{n}{2}\rceil-1$ and this value is achieved when $n$ is odd.
 \end{enumerate}
\end{proof}

In the proof of Lemma \ref{lbl-paths} it is established for a particular class of graphs, namely, finite paths, that (a), there is a labelling process that generates the minimum (maximum) number of valid paths and (b) any change to this process increases (decreases) the number of valid paths. In all of the following lemmas, property (a) is established for various other graph classes. However, for the sake of brevity, property (b) is not included as its proof in each case is similar to that for Lemma \ref{lbl-paths}.

% __________________ T R E E S
\begin{lemma}\label{lbl-trees}
If $T_n$ is a tree of finite order $n$ (a connected, acyclic $n$-graph) then ($i$) $\phi_{\min}(T_n) = 0$ and ($ii$) $\phi_{\max}(T_n) \leqslant \binom{n-1}{2}$.
\end{lemma}
\begin{proof}
\phantom{x}
\begin{enumerate}[leftmargin=*,label=(\emph{\roman*})]
 \item\label{lbl-tree-i} Suppose that an algorithm is applied to $T_n$, which iteratively
 identifies an unnumbered vertex of least degree, numbers it as the lowest of the remaining numbers and then removes it from $T_n$, along with its incident edge. Since $T_n$ is a tree, at the $i$-th step of the algorithm, for $i = 1,\dots,n$; the current subgraph of $T_n$, say $T_n(i)$, has a leaf, and such a vertex of $T_n(i)$ is numbered as the smallest possible available number. Consequently, any subsequent vertex will be numbered with a larger number. When all the vertices of the original tree $T_n$ are numbered, it will not have any valid path and thus, $\phi_{\min}(T_n) = 0$.
 \item The upper bound on $\phi_{\max}(T_n)$ can be achieved exactly  when $T_n$ is a star by iteratively numbering the vertices of $T_n$ in order of their degree, highest first. At each step, the unnumbered vertex of greatest degree is numbered with the least of the remaining numbers and is then removed from $T_n$ along with its incident edge. In this case the hub of $T_n$ is numbered as $1$, the other vertices are numbered arbitrarily and the result follows. 
\end{enumerate}
\end{proof}

% _________________ C Y C L E S
\begin{lemma}\label{lbl-cycles}
If $C_n$ is a cycle graph of finite order $n$ (a chordless, closed $n$-path) then ($i$) $\phi_{\min}(C_n) = 1$ and ($ii$) $\phi_{\max}(C_n) = \left\lceil\frac{n-1}{2}\right\rceil$.
\end{lemma}
\begin{proof}
Let $C_n$ be a cycle graph of finite order $n$. Suppose an arbitrary vertex of $C_n$, $u$ say, is numbered as 1. Since $C_n$ is a cycle, no matter how the remaining vertices of $C_n$ are numbered, the 2-path having $u$ as its middle vertex will be valid. Suppose $u$ and its two incident edges are removed from $C_n$ giving the path graph $C_n - u$.
\begin{enumerate}[leftmargin=*,label=(\emph{\roman*})]
\item\label{proof-cycle-i} By Lemma \ref{lbl-paths}-\ref{lbl-path-i}, it is possible to find a numbering of $C_n - u$ such that there are no further valid paths in $C_n$. Thus, $\phi_{\min}(C_n)=1$.
\item\label{proof-cycle-ii} By Lemma \ref{lbl-paths}-\ref{lbl-path-ii}, the maximum possible quantity of valid paths in the path graph $C_n-u$ is $\phi_{\max}(C_n - u)=\lceil\frac{n-1}{2}\rceil - 1$. Thus, $\phi_{\max}(C_n)=1+\phi_{\max}(C_n-u)=1+\lceil\frac{n-1}{2}\rceil - 1=\lceil\frac{n-1}{2}\rceil$.
\end{enumerate}
\end{proof}

\begin{helpful}\label{lbl-help2}
By Lemma \ref{lbl-cycles}, each 3-cycle in a numbered graph $G$ induces exactly one valid path in $G$. \end{helpful}

% ________________ W H E E L S
\begin{lemma}\label{lbl-wheel}
If $W_n$ is a wheel graph  of finite order $n$ (the Cartesian product of a solitary vertex and a cycle graph, i.e., $K_1\times C_{n-1}$) then ($i$) $\phi_{\min}(W_{n})=n$ and ($ii$) $\phi_{\max}(W_{n})=\binom{n-1}{2}+\left\lceil\frac{n-1}{2}\right\rceil$. 
\end{lemma}
\begin{proof}
Let $W_n$ be a wheel graph 
of finite order $n$.
\begin{enumerate}[leftmargin=*,label=(\emph{\roman*})]
\item By Lemma \ref{lbl-cycles}-\ref{proof-cycle-i} it is possible to find a numbering of the rim of $W_n$ such that the whole rim contains only one valid path. If the hub of $W_n$ is numbered as $n$, by \ref{lbl-help2}, the $(n-1)$ 3-cycles of $W_n$ each have exactly one valid path. This numbering yields exactly $n$ valid paths.

Suppose the hub had received a number $1\leqslant k <n$. Lemma \ref{lbl-cycles}-\ref{proof-cycle-i} still ensures that the rim could have been numbered with $\{1, \dots, k-1, k + 1, \dots, n\} $ to induce only one valid path. If $k=n-1$, the hub does not induce any valid paths. On the other hand,  if $1\leqslant k<n-1$, there are $\binom{n-k}{2}$ 2-paths with the hub as middle vertex.

Therefore, the proposed numbering is minimal and the result follows.

\item Suppose that the hub of $W_n$ is numbered as $1$. Then, no matter how the rim vertices of $W_n$ are numbered, there will be $\binom{n-1}{2}$ valid paths with the hub as their middle vertex and, following item \ref{proof-cycle-ii} of the proof of Lemma \ref{lbl-cycles}, there will be a further  $\lceil\frac{n-1}{2}\rceil$ valid paths around the rim. This maximal numbering achieves the stated result.
\end{enumerate}
\end{proof}

\begin{figure}[h!]
 \centering
 \subfloat[A maximal numbering with $\phi_{\max}(K_{2,4}) = 12$.]{
  \begin{tikzpicture}[
   no/.style = {draw,circle,inner sep=1pt,minimum size=7pt,fill=black},
   transform shape,
   scale = 0.8
  ]
   \useasboundingbox (-1.3,-0.2) rectangle (3.3,3.2);
  \node [no,label=right:{6}]  (a) at (2.0,0.0) {};
  \node [no,label=right:{5}]  (b) at (2.0,1.0) {};
  \node [no,label=right:{4}]  (c) at (2.0,2.0) {};
  \node [no,label=right:{3}]  (d) at (2.0,3.0) {};
  \node [no,label=left:{2}] (e) at (0.0,1.0) {};
  \node [no,label=left:{1}] (f) at (0.0,2.0) {};
   \foreach \i in {e,f}
    \foreach \j in {a,b,c,d}
     \draw (\i) -- (\j);
  \end{tikzpicture}
  \label{sfig:c-ex-5-1}
 }
 \qquad\qquad
 \subfloat[A minimal numbering with $\phi_{\max}(K_{2,4}) = 3$]{
  \begin{tikzpicture}[
   no/.style = {draw,circle,inner sep=1pt,minimum size=7pt,fill=black},
   transform shape,
   scale = 0.8
  ]
   \useasboundingbox (-1.3,-0.2) rectangle (3.3,3.2);
   \node [no,label=right:{6}]  (a) at (2.0,0.0) {};
   \node [no,label=right:{4}]  (b) at (2.0,1.0) {};
   \node [no,label=right:{2}]  (c) at (2.0,2.0) {};
   \node [no,label=right:{1}]  (d) at (2.0,3.0) {};
   \node [no,label=left:{5}] (e) at (0.0,1.0) {};
   \node [no,label=left:{3}] (f) at (0.0,2.0) {};
   \foreach \i in {e,f}
    \foreach \j in {a,b,c,d}
     \draw (\i) -- (\j);
  \end{tikzpicture}
  \label{sfig:c-ex-5-2}
 }
 \caption[]{Maximal and minimal numberings of $K_{2,4}$.}%
 \label{fig::c-ex-5}
\end{figure}
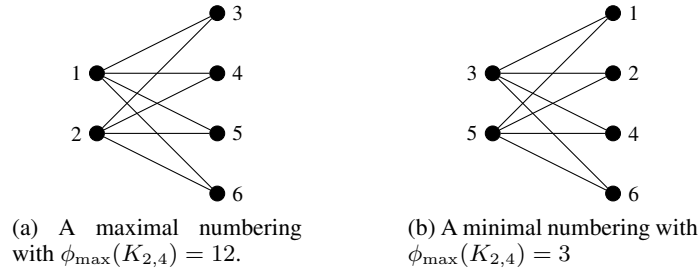

% _____ C O M P L E T E   B I P A R T I T E
\begin{lemma}\label{lbl-complete-bipartite}
If $K_{p,q}$ is the complete bipartite graph with a vertex set partition of finite orders $p$ and $q$, where $1\leqslant p \leqslant q$, then:
\begin{enumerate}[leftmargin=*,label=\emph{(}\roman*\emph{)},ref={(\emph{\roman*})},itemsep=3pt]
 \item\label{c-i} $\phi_{\max}(K_{p,q}) = p\binom{q}{2}=\frac{1}{2}p(q^2-q)$;
 \item\label{c-ii}  $\phi_{\min}(K_{p,q})=\frac{1}{6}(3q-p-1)(p^2 - p)$.
\end{enumerate}
\end{lemma}
\begin{proof}
 Let $K_{p,q}$ be the complete bipartite graph with a vertex set partition $V_p$ and $V_q$ of finite orders $p$ and $q$, respectively, where $1 \leqslant p \leqslant q$. Result \ref{c-i} follows by constructing a maximal numbering by labelling the vertices in $V_p$ with the numbers $1,2, \dots,p$; and the vertices in $V_q$ with the numbers $p+1,p+2\dots,p+q$; as illustrated in Figure  \ref{sfig:c-ex-5-1}. Result \ref{c-ii} follows by constructing a minimal numbering as follows. If $q$ is even, label the vertices in $V_p$ with the $p$ highest odd numbers $\mathcal{O}=\{o_1,o_2,\dots,o_p\}\subseteq \{1,3, \dots,2p-1\}$ and the vertices in $V_q$ with the numbers in $\{1,2,\dots,p+q\}\setminus\mathcal{O}$. 
If $q$ is odd, label the vertices in $V_p$ with the $p$ highest even numbers $\mathcal{E}=\{e_1,e_2,\dots,e_p\}\subseteq \{2,4, \dots,2p\}$ and the vertices in $V_q$ with the numbers in $\{1,2,\dots,p+q\}\setminus\mathcal{E}$.
Taking each vertex of one partition as the middle vertex of a possible valid path and combining it with each pair of higher-labelled vertices of the other partition, allows the following development:
 \begin{align*}
    \phi_{\min}(K_{p,q}) & =
     \begin{cases}
      0,   & \text{if } p=1 \text{ and } q\geqslant 1;\\
      q-1, & \text{if } p=2 \text{ and } q\geqslant 2;\\
      3q-4, & \text{if } p=3 \text{ and } q\geqslant 3;\\
      6q-10, & \text{if } p=4 \text{ and } q\geqslant 4;\\
      10q-20, & \text{if } p=5 \text{ and } q\geqslant 5;\\
      15q-35, & \text{if } p=6 \text{ and } q\geqslant 6;\\
    %   21q-56, & \text{if } p=7 \text{ and } q\geqslant 7;\\
    %   \text{and so forth for} & p\geqslant 8.\\
    %   \text{and so forth for} & p\geqslant 7 \text{ and } q\geqslant p.\\
    \end{cases}\shortintertext{and so forth for $p\geqslant 7$ and $q\geqslant p$. This result can be generalised as}
    \phi_{\min}(K_{p,q}) & =\frac{1}{2}(p^2 - p)q - \frac{1}{6}p(p^2-1),\shortintertext{for $p\geqslant 1$ and $q\geqslant p$. Some simple algebraic manipulations, gives}  
    \phi_{\min}(K_{p,q}) & = \frac{1}{6}(3q-p-1)(p^2 - p),\shortintertext{for $p\geqslant 1$ and $q\geqslant p$.}
 \end{align*}
\end{proof}

The first few cases for Result \ref{c-ii} can be proved as follows. When $p=1$, if the number $q$ is assigned to $V_1$, then $\phi_{\min}(K_{1,q}) = 0$. When $p=2$, if the two highest odd numbers in $\{1,2,\dots,q+2\}$ are assigned to $V_1$, then $\phi_{\min}(K_{2,q}) = q-1$. When $p=3$ and $q=4$, if the numbers in $\{3,5,7\}$ are assigned to $V_1$, then $\phi_{\min}(K_{3,4}) = 8$. When $p=3$ and $q\geqslant 5$, if the three highest odd numbers in $\{1,2,\dots,q+3\}$ are assigned to $V_1$, then $\phi_{\min}(K_{3,q}) = \phi_{\min}(K_{3,q-1}) + 3$. Result \ref{c-ii} is illustrated in Figure \ref{sfig:c-ex-5-2}. ($\phi_{\min}(K_{3,q}) = (5,8,11,14,17,20,23,\dots)$, for $q=1,2,3,4,5,6,7,\dots$).

Case $p = q$ of Result \ref{c-ii} follows by counting the number of valid paths obtained by taking each vertex of a partition as the middle vertex and combining it with each pair of vertices with higher labels from the other partition. Suppose that the labels of the vertices in $V_p$ are odd and those associated to the vertices in $V_q$ are even. Thus, The smallest of the labels in $V_p$ is smaller than all the labels in $V_q$, the second smallest label in $V_p$ is greater than just one of the labels in $V_q$, the third smallest label is greater than two of the labels in $V_q$, and so on. Therefore, the number of valid paths is given by:
\begin{align*}
\phi_{\min}(K_{q,q}) & = \binom{q}{2} + 2\binom{q-1}{2} + 2\binom{q-2}{2} +\cdots + 2\binom{2}{2}\\
 & = \binom{q}{2} + 2\sum_{i=1}^{q-2}\binom{q-i}{2}\\
 & = \tfrac{1}{6}(q-1)(2q^2-q) = \tfrac{1}{6}(2q-1)(q^2 - q)\\
 & = (0,1,5,14,30,55,\dots),\textbf{ for }  q=1,2,3,4,5,6,\dots
\end{align*}

%\begin{corollary}
%If $K_{p,p}$ is the complete bipartite graph with finite order $p$, then $\phi_{\min}(K_{p,p}) = \phi_{\max}(K_{p,p})$.
%\end{corollary}

\begin{helpful}\label{lbl-help1}
Every vertex of $K_n$, the complete graph of finite order $n$, is the middle vertex of $\binom{n-1}{2}$ 2-paths.
\end{helpful}

% _________ C O M P L E T E   G R A P H
\begin{lemma}\label{lbl-complete-graph}
If $K_n$ is the complete graph of finite order $n$, any numbering of $K_n$ induces 
$\tfrac{1}{6}(n^3 - 3n^2 +2n)$ valid 2-paths. 
\end{lemma}
\begin{proof}
Suppose $\pi$ is a numbering of $K_n$. Consider the sequence $S_{\pi} = \langle u_1, u_2,\dots,\allowbreak u_n\rangle$, of the vertices of $K_n$ assembled in order of their numbers assigned by $\pi$, lowest first, i.e., $\pi(u_i) = i$, $i = 1,\dots,n$. By \ref{lbl-help1}, as $f(u_1) = 1$, the minimum vertex number, $u_1$ is the middle vertex of $\binom{n-1}{2}$ valid paths. No other vertex can be the middle vertex of any valid path with $u_1$ as one of the endpoints. Hence, without danger of miscounting, we can consider $u_2$ in $K_n - u_1$, which is a complete graph on $n-1$ vertices. By an analogous argument to that just made, $u_2$ is the middle vertex of $\binom{n-2}{2}$ valid paths in the remaining $K_{n-2}$. Continuing in this way, it follows easily that $u_i$ is the middle vertex of $\binom{n-i}{2}$ valid paths, $i = 1,\dots,n$. Hence the total number of valid paths is $\sum^{n-1}_{i = 2}\binom{i}{2} = \tfrac{1}{6}(n^3 - 3n^2 + 2n)$.
\end{proof}

As the result in Lemma \ref{lbl-complete-graph} is numbering-independent, Corollary \ref{cor:Kn} below  follows immediately.
\begin{corollary}\label{cor:Kn}
If $K_n$ is the complete graph of finite order $n$ and $\pi$ is any numbering of $K_n$, the minimum quantity of valid paths equals the maximum quantity of valid paths, i.e.,
\begin{equation}\label{eq:cor-complete}
\phi_\pi(K_n) = \phi_{\min}(K_n) = \phi_{\max}(K_n) = \tfrac{1}{6}(n^3 - 3n^2 +2n).
\end{equation}
\end{corollary}

Regarding \ref{lbl-help2}, planar graphs of five or more vertices often contain two types of 3-cycles -- faces and separating triangles, which aids in the analysis of numberings of graphs in this class. However, these and many other graphs have additional valid paths that are not members of 3-cycles. Of course, there are many graphs that are 3-cycle-free, including trees, $n$-cycles, bipartite graphs and the connected strongly regular graphs (the 5-cycle, all complete bipartite graphs and the Petersen, Clebsch, Hoffmann-Singleton, Gerwirtz and M22 graphs).

A planar graph $G$ is said to be a maximal planar graph (MPG) if the addition of any edge to $G$ results in a nonplanar graph. Suppose that $G$ is an MPG. Let $\tau(G)$ be the number of triangles (3-cycles) in $G$. A separating triangle of $G$ is a triangle of $G$ that is not a face of $G$. Let $\sigma(G)$ be the number of separating triangles in $G$. The operation of inserting exactly one new vertex in a face of $G$ and joining the new vertex to the three vertices of the face is termed \emph{dimpling} the face \cite{Foulds-Robinson-1979}.

\begin{theorem}[\cite{Hakimi-Schmeichel-1979}]\label{teo:hakimi}
If $G$ is an MPG with more than $5$ vertices ($n > 5$) then $2n - 4 \leqslant \tau(G) \leqslant 3n - 8$.
\end{theorem}

\begin{proof}
The lower bound is obvious as $G$ is an MPG and thus has $2n - 4$ faces, all of which are triangles. (Incidentally, a necessary and sufficient condition for the lower bound to be tight (that is, $\tau(G) = 2n - 4)$ is that $G$ is 4-connected.) The upper bound is achieved if and only if $G$ is obtained from $K_4$ by iteratively dimpling faces of $K_4$ and the resulting MPG's. This produces $n-4$ separating triangles.
\end{proof}

\begin{corollary}[\cite{Helden-2007}]\label{cor:helden}
If $G$ is an MPG with more than 5 vertices then $0 \leqslant \sigma(G) \leqslant n - 4$.
\end{corollary}

\begin{proof}
From Theorem \ref{teo:hakimi}, $\tau(G) \leqslant 3n - 8$. Furthermore, as $G$ is an MPG it has $2n-4$ faces. By definition, any separating triangle is not a face. So $\sigma(G)$ can be at most $n - 4$.
\end{proof}

By its very nature, a sequence of face dimpling operations can be used to construct only a strict subset of all MPGs, called Apollonian graphs \cite{Andrade-Herrmann-Andrade-Silva-2005}. For example, the octahedral graph (shown in Figure 7(a)) does not have a vertex of degree 3 and thus cannot be generated by a sequence of only face dimpling operations.

% _____________ A P O L L O N I A N
\begin{lemma}\label{lbl-apollonian}
If $A_n$ is an Apollonian graph of finite order $n$, then 
\begin{enumerate}[leftmargin=*,label=\emph{(}\roman*\emph{)},ref={(\emph{\roman*})},itemsep=3pt]
 \item\label{a-i} $\phi_{\min}(A_n) = 3n-8$, and
 \item\label{a-ii} 
$\phi_{\max}(A_n) = 
   \begin{cases}
(a)\ 4\left(\frac{n}{2}-1\right)^2+\frac{n}{2}-2 =  n^2-\frac{7}{2}n+2, &\text{if } n\geqslant 4 \text{ and even;}\\
(b)\ 4\left(\left\lfloor\frac{n}{2}\right\rfloor\right)^2 - 3\left\lfloor\frac{n}{2}\right\rfloor, &\text{if } n\geqslant 3 \text{ and odd.}
\end{cases}$
\end{enumerate}
\end{lemma}

It follows from \ref{a-ii} that $\phi_{\max}(A_n)= 1, 4, 10, 17, 27,38, 52,\dots$, for $n=3,4,5,\allowbreak6,7,8,9,\dots$

%\les{\begin{helpful}\label{lbl-help3}
%Note that the relations $(ii)$ in the statement of Lemma \ref{lbl-apollonian} always hold when $A_n$ is an Apollonian graph of the form $P_2\times P_{n-2}$, where $P_2$ and $P_{n-2}$ are distinct paths of length 2 and $n-2$, respectively. Thus, such a graph is maximal in the sense that any numbering of any of the other Apollonian graphs of order $n$ will produce a quantity of valid paths that is no greater than either $(n^2-\tfrac{7}{2}n+2)$ when $n$ is even, or $(4\left(\left\lfloor\frac{n}{2}\right\rfloor\right)^2 - 3\left\lfloor\frac{n}{2}\right\rfloor)$ when $n$ is odd.
%\end{helpful}

\begin{proof}
Let $A_n$ be an Apollonian graph of finite order $n$ with a numbering $\pi$. Clearly, $\pi$ induces exactly one valid path for each of the 3-cycles of $A_n$.
\begin{enumerate}[leftmargin=*,label=(\emph{\roman*})]
\item Suppose $\pi$ corresponds to an application of an algorithm which progressively numbers the vertices of $A_n$ in the reverse order to that being used to construct $A_n$ by face dimpling. In this case, it is straightforward to show by mathematical induction (given in the next paragraph) that there are no further valid paths in $A_n$ beyond those induced by its triangles. As $A_n$ is an Apollonian graph, the result follows from the last observation in the proof of Theorem \ref{teo:hakimi}.

\medskip\emph{Proof by induction:}\\
Let $P(n)$ be the mathematical statement $\phi_{\min}(A_n) = 3n-8$.
\begin{description}[noitemsep, topsep=0pt,font=\normalfont]%, before={\vspace*{-\baselineskip}}]
\item [\emph{Base case:}] When $n = 4$ we have $A_4 = K_4$. By Lemma \ref{lbl-complete-graph}, $\phi_{\min}(A_n) = \frac{4^3 -  3{\cdot}4^2 + 2{\cdot}4}{6} = 4 = 3{\cdot}4 -8$. So $P(4)$ is correct.
\item [\emph{Inductive hypothesis:}] Assume that $P(k)$ is correct for some integer $k \geqslant 4$. This implies that $\phi_{\min}(A_k) = 3k-8$.
\item [\emph{Inductive step:}] We now show that $\phi_{\min}(A_{k+1}) = 3(k+1) - 8 = 3k - 5$. As $A_{k+1}$ is an Apollonian graph it can be constructed by face dimpling. Let the last vertex (which is of degree $3$) to be dimpled into $A_{k+1}$ be denoted by $u$. Consider $A_{k+1} - {u}$. By the inductive hypothesis
\begin{equation}\label{eq-proof}
\phi_{\min}(A_{k+1} -{u}) = 3k - 8
\end{equation} and there exists a numbering with $\{1,\dots,k\}$. Renumber the vertices by increasing each number by unity. This preserves the correctness of \eqref{eq-proof}. Return $u$ to the face it was in to re-establish $A_{k+1}$ and number $u$ as unity. Then $\phi_{\min}(A_{k+1}) = 3k - 8 +3 = 3(k+1) - 8 = 3k - 5$. This completes the proof of $(i)$.
\end{description}

\item The first few cases follow by exhaustive analysis, including the only Apollonian graphs with $n = 3,4,5,6,7$, being respectively: for $n=3$, the triangle (3-cycle); for $n=4$, $K_4$ (the complete graph with four vertices); for $n=5$, the Johnson solid skeleton $12$ (Figure \ref{sfig:johnson-12}); for $n=6$, the Hexahedral graph 5, (Figure \ref{sfig:hexahedral-5}); and for $n=7$, the three graphs shown in Figure \ref{fig:heptahedrals}: the heptahedral graph 15, the heptahedral graph 29 and the 2-Apollonian graph. For any given Apollonian graph of order $n$, the maximal value of $\phi_{\max}(A_{n})$ can be attained by numbering its vertices as $1,2,\dots,n$; in the order of vertex degree, starting with highest degree first. Ties can be broken arbitrarily. %For example, the graph in Figure \ref{sfig:heptahedral-15} is maximal in the sense that it can be numbered in the way just mentioned to achieve $\phi_{\max}(A_{7}) = 27$, the highest number of valid paths among all Apollonian graphs of order seven. The maximal numberings of the graphs in Figures \ref{sfig:heptahedral-29} and \ref{sfig:item-c} achieve the maximum numbers of valid paths in these graphs, being \les{$24$ and $25$, respectively.} %The following inductive proofs of the general cases of $(ii)(a)$ and $(ii)(b)$ rely on the results in \ref{lbl-help3}. Observe that this is sufficient as the graphs described there are maximal in the sense that no different graph can induce more valid paths.
\begin{figure}[!htbp]
 \centering
 \subfloat[Johnson solid skeleton 12 ($n=5$).]{
 \begin{tikzpicture}[%FIGURE 3 (A)
  no/.style = {draw,circle,inner sep=2pt,minimum size=5pt,fill=black},
%  every label/.append style={label distance=-2pt,text=red, font=\footnotesize},
  transform shape,
  scale = 0.5
 ]
  \node[no] (a) at (0.00,0.00) {};
  \node[no] (b) at (4.00,0.00) {};
  \node[no] (c) at (2.00,3.46) {};
  \node[no] (d) at (2.00,2.21) {};
  \node[no] (e) at (2.00,0.96) {};
  \draw (a) to [out=-15,in=195] (b)
        (b) to [out=105,in=-30] (c)
        (c) to [out=210,in=75]  (a);
  \foreach [count=\c from 0] \d in {e,d}
  {
   \draw (a) to [out=30+25*\c,in=185.25+6*\c] (\d);
   \draw (b) to [out=150-25*\c,in=-5.25-6*\c] (\d);
  }
  \draw (c) -- (e);
 \end{tikzpicture}
  \label{sfig:johnson-12}
 }
 \qquad\qquad
 \subfloat[Hexahedral graph 5 ($n=6$).]{
 \begin{tikzpicture}[%FIGURE 3 (B)
  no/.style = {draw,circle,inner sep=2pt,minimum size=5pt,fill=black},
%  every label/.append style={label distance=-2pt,text=red, font=\footnotesize},
  transform shape,
  scale = 0.5
 ]
  \node[no] (a) at (0.00,0.00) {};
  \node[no] (b) at (4.00,0.00) {};
  \node[no] (c) at (2.00,3.46) {};
  \node[no] (d) at (2.00,2.52) {};
  \node[no] (e) at (2.00,1.68) {};
  \node[no] (f) at (2.00,0.74) {};
  \draw (a) to [out=-15,in=195] (b)
        (b) to [out=105,in=-30] (c)
        (c) to [out=210,in=75]  (a);
  \foreach [count=\c from 0] \d in {f,e,d}
  {
   \draw (a) to [out=20+20*\c,in=185.25+6*\c] (\d);
   \draw (b) to [out=160-20*\c,in=-5.25-6*\c] (\d);
  }
  \draw (c) -- (f);
 \end{tikzpicture}
  \label{sfig:hexahedral-5}
 }
 \caption{The only Apollonian graphs with $n = 5$ or $6$.}
 \label{fig:item-a-b}
\end{figure}
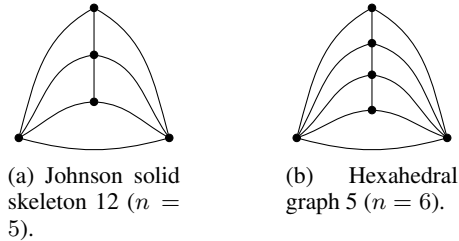
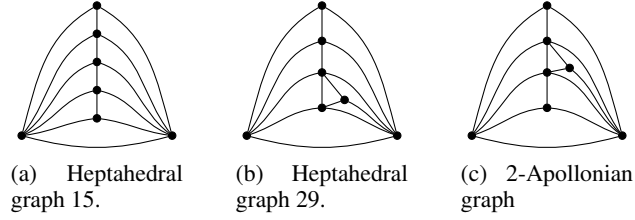
\begin{figure}[!htbp]
 \centering
 \subfloat[Heptahedral graph 15.]{
 \begin{tikzpicture}[%FIGURE 4 (A)
  no/.style = {draw,circle,inner sep=2pt,minimum size=4pt,fill=black},
%  every label/.append style={label distance=-1pt,text=red, font=\footnotesize},
  transform shape,
  scale = 0.5
 ]
  \node[no] (a) at (0.00,0.00) {};
  \node[no] (b) at (4.00,0.00) {};
  \node[no] (c) at (2.00,3.46) {};
  \node[no] (d) at (2.00,2.71) {};
  \node[no] (e) at (2.00,1.96) {};
  \node[no] (f) at (2.00,1.21) {};
  \node[no] (g) at (2.00,0.46) {};
  \draw (a) to [out=-15,in=195] (b)
        (b) to [out=105,in=-30] (c)
        (c) to [out=210,in=75]  (a);
  \foreach [count=\c from 0] \d in {g,f,e,d}
  {
   \draw (a) to [out=7.5+15*\c,in=187.5+6*\c]  (\d);
   \draw (b) to [out=172.5-15*\c,in=-7.5-6*\c] (\d);
  }
  \draw (c) -- (g);
 \end{tikzpicture}
  \label{sfig:heptahedral-15}
 }\qquad
 \subfloat[Heptahedral graph 29.]{
 \begin{tikzpicture}[%FIGURE 4 (B)
  no/.style = {draw,circle,inner sep=2pt,minimum size=4pt,fill=black},
%  every label/.append style={label distance=-1pt,text=red, font=\footnotesize},
  transform shape,
  scale = 0.5
 ]
  \node[no] (a) at (0.00,0.00) {};
  \node[no] (b) at (4.00,0.00) {};
  \node[no] (c) at (2.00,3.46) {};
  \node[no] (d) at (2.00,2.52) {};
  \node[no] (e) at (2.00,1.68) {};
  \node[no] (f) at (2.00,0.74) {};
  \node[no] (g) at (2.60,0.95) {};
  \draw (a) to [out=-15,in=195] (b)
        (b) to [out=105,in=-30] (c)
        (c) to [out=210,in=75]  (a);
  \foreach [count=\c from 0] \d in {f,e,d}
  {
   \draw (a) to [out=18.7+18*\c,in=185.25+6*\c]  (\d);
   \draw (b) to [out=161.3-18*\c,in=-5.25-6*\c] (\d);
  }
  \draw (c) -- (f) -- (g) -- (e) (b) to [out=154,in=-25] (g);

 \end{tikzpicture}
  \label{sfig:heptahedral-29}
 }\qquad
 \subfloat[2-Apollonian graph]{
 \begin{tikzpicture}[%FIGURE 4 (C)
  no/.style = {draw,circle,inner sep=2pt,minimum size=4pt,fill=black},
%  every label/.append style={label distance=-1pt,text=red, font=\footnotesize},
  transform shape,
  scale = 0.5
 ]
  \node[no] (a) at (0.00,0.00) {};
  \node[no] (b) at (4.00,0.00) {};
  \node[no] (c) at (2.00,3.46) {};
  \node[no] (d) at (2.00,2.52) {};
  \node[no] (e) at (2.00,1.68) {};
  \node[no] (f) at (2.00,0.74) {};
  \node[no] (g) at (2.60,1.80) {};
  \draw (a) to [out=-15,in=195] (b)
        (b) to [out=105,in=-30] (c)
        (c) to [out=210,in=75]  (a);
  \foreach [count=\c from 0] \d in {f,e,d}
  {
   \draw (a) to [out=18.7+18*\c,in=185.25+6*\c]  (\d);
   \draw (b) to [out=161.3-18*\c,in=-5.25-6*\c] (\d);
  }
  \draw (c) -- (f) (d) -- (g) -- (e) (b) to [out=135,in=-40] (g);
 \end{tikzpicture}
  \label{sfig:item-c}
 }
 \caption{The Apollonian graphs with $n = 7$.}
 \label{fig:heptahedrals}
\end{figure}

Firstly, we consider the special case where $A_n$ is an Apollonian graph of the form $P_2\times P_{n-2}$, where $P_2$ and $P_{n-2}$ are distinct paths of length 2 (i.e., with 2 vertices) and $(n-2)$ (i.e., with $(n-2)$ vertices), respectively. The 2-path is denoted by $\langle v_1, v_2\rangle$ and is termed the \emph{base}, and the $(n-2)$-path is denoted by $\langle v_3, v_4,\dots v_{n}\rangle$ and is termed the \emph{spine}. An example of such a graph is illustrated in Figure \ref{sfig:hexahedral-5}. It is straightforward to establish that any maximal numbering of $A_n$ is such that the base is numbered with $1$ and $2$ (in either order). In this case, all triples consisting of the vertex of the base that is numbered 1, and any other two vertices of $A_n$ (including the vertex numbered 2) induce $\binom{n-1}{2}$  %$C\rlap{\textsuperscript{n-1}}\textsubscript{2}$   
valid paths. Similarly,  all triples consisting of a pair of members of the spine and the vertex of the base that is numbered 2, induce 
$\binom{n-2}{2}$ 
% $C\rlap{\textsuperscript{n-2}}\textsubscript{2}$   
valid paths. Therefore, there are at least
$\binom{n-1}{2} + \binom{n-2}{2}$  
% $C\rlap{\textsuperscript{n-1}}\textsubscript{2}$ + $C\rlap{\textsuperscript{n-2}}\textsubscript{2}$ 
valid paths. The question is what is the maximum number of additional valid paths that can be induced by triples of vertices solely within the spine (which has $(n-2)$ vertices)? In other words, how should the vertices of the spine be numbered so as to maximise the total number of valid paths? By Lemma \ref{lbl-paths}, if $n$ is even there is a maximum of $\lceil\frac{n-2}{2}\rceil-1$ additional valid paths, and if $n$ is odd there is a maximum of $\lfloor\frac{n-2}{2}\rfloor$ additional paths. These results imply that:
\begin{align*}
\phi_{\max}(P_2\times P_{n-2}) & = (n-1)(n-2)/2+(n-2)(n-3)/2 + \lceil\frac{n-2}{2}\rceil-1\\ & = n^2-\tfrac{7}{2}n+2, \text{ for } n\geqslant 4 \text{ and even,}\shortintertext{and also}
\phi_{\max}(P_2\times P_{n-2}) & = (n-1)(n-2)/2+(n-2)(n-3)/2 + \lfloor\tfrac{n-2}{2}\rfloor\\ & = n^2-4n+3+\lfloor\frac{n}{2}\rfloor\\ & = 4\left(\left\lfloor\tfrac{n}{2}\right\rfloor\right)^2-3\left\lfloor\tfrac{n}{2}\right\rfloor,\text{ for } n\geqslant 3 \text{ and odd.}\\
\end{align*}

Examples of maximal numberings of such $P_2\times P_{n-2}$ graphs are given in Figure \ref{fig:P2Pn-2}. The numbering of the graph in Figure \ref{sfig:p2-p6}, where $n=8$ (even), induces 38 valid paths, the maximum among all graphs of order 8. The numbering of the graph in Figure \ref{sfig:p2-p7}, where $n=9$ (odd), induces 52 valid paths, the maximum among all graphs of order 9.
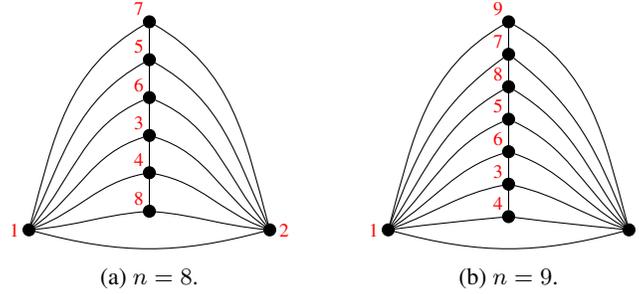
\begin{figure}[!htbp]
 \centering
 \subfloat[$n=8$.]{
 \begin{tikzpicture}[%FIGURE 5 (A)
  no/.style = {draw,circle,inner sep=2pt,minimum size=4pt,fill=black},
  every label/.append style={label distance=-4pt,text=red, font=\footnotesize},
  transform shape,
  scale = 0.8
 ]
  \node[no,label={[label distance=-2pt]180:1}] (a) at (0.00,0.00) {};
  \node[no,label={[label distance=-2pt]000:2}] (b) at (4.00,0.00) {};
  \node[no,label=150:{7}] (c) at (2.00,3.46) {};
  \node[no,label=150:{5}] (d) at (2.00,2.83) {};
  \node[no,label=150:{6}] (e) at (2.00,2.20) {};
  \node[no,label=150:{3}] (f) at (2.00,1.57) {};
  \node[no,label=150:{4}] (g) at (2.00,0.95) {};
  \node[no,label=150:{8}] (h) at (2.00,0.31) {};
  \draw (a) to [out=-15,in=195] (b)
        (b) to [out=105,in=-30] (c)
        (c) to [out=210,in=75]  (a);
  \foreach [count=\c from 0] \d in {h,g,f,e,d}
  {
   \draw (a) to [out=5.25+15*\c,in=185.25+6*\c]  (\d);
   \draw (b) to [out=174.75-15*\c,in=-5.25-6*\c] (\d);
  }
  \draw (c) -- (h);
 \end{tikzpicture}
  \label{sfig:p2-p6}
 }
 \qquad
 \subfloat[$n=9$.]{
 \begin{tikzpicture}[%FIGURE 5 (B)
  no/.style = {draw,circle,inner sep=2pt,minimum size=4pt,fill=black},
  every label/.append style={label distance=-4pt,text=red, font=\footnotesize},
  transform shape,
  scale = 0.8
 ]
  \node[no,label={[label distance=-2pt]180:1}] (a) at (0.00,0.00) {};
  \node[no,label={[label distance=-2pt]000:2}] (b) at (4.00,0.00) {};
  \node[no,label=150:{9}] (c) at (2.00,3.46) {};
  \node[no,label=150:{7}] (d) at (2.00,2.92) {};
  \node[no,label=150:{8}] (e) at (2.00,2.38) {};
  \node[no,label=150:{5}] (f) at (2.00,1.84) {};
  \node[no,label=150:{6}] (g) at (2.00,1.30) {};
  \node[no,label=150:{3}] (h) at (2.00,0.76) {};
  \node[no,label=150:{4}] (i) at (2.00,0.22) {};
  \draw (a) to [out=-15,in=195] (b)
        (b) to [out=105,in=-30] (c)
        (c) to [out=210,in=75]  (a);
  \foreach [count=\c from 0] \d in {i,h,g,f,e,d}
  {
   \draw (a) to [out=7.5+12.5*\c,in=185.5+6*\c]  (\d);
   \draw (b) to [out=174.5-12.5*\c,in=-7.25-6*\c] (\d);
  }
  \draw (c) -- (i);
 \end{tikzpicture}
  \label{sfig:p2-p7}
 }
 \caption{Examples of some maximal numberings of $P_2\times P_{n-2}$ graphs.}
 \label{fig:P2Pn-2}
\end{figure}

Secondly, we consider the rest of the cases, i.e. the set of Apollonian graphs of order $n$ that are not of the form $P_2\times P_{n-2}$. Let $Q_{2,n-2}$ denote this set. Note that if any such graph in $Q_{2,n-2}$ has a $P_2$ path and a $P_{n-2}$ path, then at least one of the vertices in its $P_{n-2}$ path is not adjacent to one of the vertices in its $P_2$ path. The extreme case (in the sense that graphs in $Q_{2,n-2}$ that can be numbered so as to induce the most valid paths among all graphs in $Q_{2,n-2}$) occurs for graphs that have both a $P_2$ and a $P_{n-2}$ and there is exactly one adjacency missing between the two paths. Let $H_n = G_{2,P_{n-2}}$ be such a graph. We show that 
\begin{align}
\phi_{\max}(H_n) & \leqslant n^2-\tfrac{7}{2}n+2,\text{ for } n\geqslant 4 \text{ and even.}\label{Qna} \\\shortintertext{and}
\phi_{\max}(H_n) & \leqslant 4\left(\left\lfloor\tfrac{n}{2}\right\rfloor\right)^2-3\left\lfloor\tfrac{n}{2}\right\rfloor,\text{ for }  n\geqslant 3 \text{ and odd.}\label{Qnb}
\end{align}

\paragraph{Proof of \eqref{Qna}:} The case when $n=4$ holds by inspection. Now consider $G_{2\times P_{n-2}}$, where $n\  (\geqslant 6)$ is even. Suppose that vertex $v_k$, $3\leqslant k \leqslant n$, is the unique vertex in $P_{n-2}$ that is not adjacent to one of the vertices of the base, $v_2$ say. The extreme case occurs when $k=n$, i.e. $v_n$ is an endpoint of $P_{n-2}$ (and is thus adjacent to $v_{n-1}$) but is also adjacent to only $v_1$ and $v_{n-2}$. The maximal numbering occurs when the vertices are numbered as they are for $P_2\times P_{n-2}$. (Even though $v_n$ has been relocated from face $\langle v_1,v_2,v_n\rangle$ to face $\langle v_1,v_{n-1},v_{n-2}\rangle$, it is assigned the same number.) Recall that as $n$ is even, when $P_2\times P_{n-2}$ is maximally numbered,
$\binom{n-1}{2}+\binom{n-2}{2}+\lceil\frac{n-2}{2}\rceil-1$
% $(C\rlap{\textsuperscript{n-1}}\textsubscript{2}$ + $C\rlap{\textsuperscript{n-2}}\textsubscript{2}$ + $\lceil\frac{n-2}{2}\rceil-1)$
valid paths are induced. In present case, compared to $P_2\times P_{n-2}$, the $(n-3)$ paths $\langle k,2,n\langle, k=3,4,\dots,n-1$; are no longer valid. Although the relocation of $v_n$ may induce up to two new valid paths, the total number of valid paths that are invalidated always exceeds two and (\ref{Qna}) is proven.

\paragraph{Proof of \eqref{Qnb}:} Consider $G_{2\times P_{n-2}}$, where $n$ is odd. The reasoning is similar to that used to prove (\ref{Qna}) above. The case when $n=3$ holds by inspection. Now consider $G_{2\times P_{n-2}}$, where $n\ (\geqslant 5)$ is odd. Suppose that vertex $v_k$, $3\leqslant k \leqslant n$, is the unique vertex in $P_{n-2}$ that is not adjacent to one of the vertices of the base, $v_2$ say. The extreme case occurs when $k=n$, i.e. $v_n$ is an endpoint of $P_{n-2}$ (and is thus adjacent to $v_{n-1}$) but is also adjacent to only $v_1$ and $v_{n-2}$. The maximal numbering occurs when the vertices are numbered as they are for $P_2\times P_{n-2}$. (Even though $v_n$ has been relocated from face $\langle v_1,v_2,v_n\rangle$ to face $\langle v_1,v_{n-1},v_{n-2}\rangle$, it is assigned the same number.) Recall that as $n$ is odd, when $P_2\times P_{n-2}$ is maximally numbered,
$\binom{n-1}{2}+\binom{n-2}{2}+\lfloor\frac{n-2}{2}\rfloor)$
% $(C\rlap{\textsuperscript{n-1}}\textsubscript{2}$ + $C\rlap{\textsuperscript{n-2}}\textsubscript{2}$ + $\lfloor\frac{n-2}{2}\rfloor)$
valid paths are induced. In present case, compared to $P_2\times P_{n-2}$, the $(n-3)$ paths $\langle k,2,n\rangle, k=3,4,\dots,n-1$; are no longer valid. Although the relocation of $v_n$ may induce up to two new valid paths, the total number of valid paths that are invalidated always exceeds two and \eqref{Qnb} is proven.\medskip

It follows from \eqref{Qna} and \eqref{Qnb} that when $A_n$ is of the form $P_2\times P_{n-2}$, it is maximal in the sense that a numbering of it can be found that induces the maximum quantity of valid paths among all numberings of all Apollonian graph of order $n$. This completes the proof.
\end{enumerate}
\end{proof}

A first example of \eqref{Qnb}, without a numbering, is illustrated in Figure \ref{fig:heptahedrals}. The Apollonian graph (of the form $P_2\times P_{n-2}$) given in Figure \ref{sfig:heptahedral-15} is maximal in the sense that it can be numbered to achieve $\phi_{\max}(A_{7}) = 27$, the highest quantity of valid paths among all Apollonian graphs of order 7. The other Apollonian graphs of order 7 (which are not of the form $P_2\times P_{n-2}$) are given in Figures \ref{sfig:heptahedral-29} and \ref{sfig:item-c}. Their maximal numberings achieve lesser quantities of valid paths namely, $24$ and $25$, respectively.

Further examples of the suboptimal numberings of such $Q_{2,n-2}$ graphs are given in Figure \ref{fig:Q2n-2}. The numbering of the graph in Figure \ref{sfig:q-2-6}, where $n=8$ (even), induces $35$ valid paths, which is less than the maximum of $\phi_{\max}(A_{8}) = 38$, illustrated in Figure \ref{sfig:p2-p6}. The numbering of the graph in Figure \ref{sfig:q-2-6}, where $n=9$ (odd), induces $46$ valid paths, which is less than the maximum of $\phi_{\max}(A_{9}) = 52$, illustrated in Figure \ref{sfig:p2-p7}.

The results in Lemma \ref{lbl-apollonian} are not valid for all MPGs. For example, when $G$ is the (nonApollonian) octahedral graph, $n = 6$, $\phi_{\min}(A_n) = 11 > 3n - 8$ and $\phi_{\max}(A_n) = 14 < 17$. This counterexample is shown in Figures \ref{sfig:octahedral-i} and \ref{sfig:octahedral-ii}.  Similar phenomena occur when $n=7$, and $G$ is either the heptahedral graph 34 or the Johnson solid skeleton 13. These counterexamples are shown in Figures  \ref{sfig:graph34-i}, \ref{sfig:johnson13-i},  \ref{sfig:graph34-ii} and \ref{sfig:johnson13-ii}.

\begin{figure}[!htbp]
 \centering
 \subfloat[$n=8$.]{
 \begin{tikzpicture}[%FIGURE 6 (A)
  no/.style = {draw,circle,inner sep=1pt,minimum size=4pt,fill=black},
  every label/.append style={label distance=-4pt,text=red, font=\footnotesize},
  transform shape,
  scale = 0.75
 ]
  \node[no,label={[label distance=-2pt]180:2}] (a) at (0.00,0.00) {};
  \node[no,label={[label distance=-2pt]000:1}] (b) at (4.00,0.00) {};
  \node[no,label=150:{7}] (c) at (2.00,3.46) {};
  \node[no,label=150:{5}] (d) at (2.00,2.71) {};
  \node[no,label=150:{6}] (e) at (2.00,1.96) {};
  \node[no,label=150:{3}] (f) at (2.00,1.21) {};
  \node[no,label=150:{4}] (g) at (2.00,0.46) {};
  \node[no,label={[label distance=-2.5pt]90:8}] (h) at (2.60,0.6) {};
  \draw (a) to [out=-15,in=195] (b)
        (b) to [out=105,in=-30] (c)
        (c) to [out=210,in=75]  (a);
  \foreach [count=\c from 0] \d in {g,f,e,d}
  {
   \draw (a) to [out=7.5+15*\c,in=187.5+6*\c]  (\d);
   \draw (b) to [out=172.5-15*\c,in=-7.5-6*\c] (\d);
  }
  \draw (c) -- (g) -- (h) -- (f) (b) to [out=167,in=-12] (h);
 \end{tikzpicture}
  \label{sfig:q-2-6}
 }
 \qquad
 \subfloat[$n=9$.]{
 \begin{tikzpicture}[%FIGURE 6 (B)
  no/.style = {draw,circle,inner sep=1pt,minimum size=4pt,fill=black},
  every label/.append style={label distance=-4pt,text=red, font=\footnotesize},
  transform shape,
  scale = 0.75
 ]
  \node[no,label={[label distance=-2pt]180:2}] (a) at (0.00,0.00) {};
  \node[no,label={[label distance=-2pt]000:1}] (b) at (4.00,0.00) {};
  \node[no,label=150:{4}] (c) at (2.00,3.46) {};
  \node[no,label=150:{3}] (d) at (2.00,2.83) {};
  \node[no,label=150:{6}] (e) at (2.00,2.20) {};
  \node[no,label=150:{5}] (f) at (2.00,1.57) {};
  \node[no,label=150:{8}] (g) at (2.00,0.94) {};
  \node[no,label=150:{7}] (h) at (2.00,0.31) {};
  \node[no,label={[label distance=-5pt]060:9}] (i) at (2.60,0.41) {};
  \draw (a) to [out=-15,in=195] (b)
        (b) to [out=105,in=-30] (c)
        (c) to [out=210,in=75]  (a);
  \foreach [count=\c from 0] \d in {h,g,f,e,d}
  {
   \draw (a) to [out=5.25+15*\c,in=185.25+6*\c]  (\d);
   \draw (b) to [out=174.75-15*\c,in=-5.25-6*\c] (\d);
  }
  \draw (c) -- (h) -- (i) -- (g) (b) to [out=167,in=-12] (i);
 \end{tikzpicture}
  \label{sfig:q-2-7}
 }
 \caption{Examples of some maximal numberings of $Q_{2,n-2}$ graphs.}
 \label{fig:Q2n-2}
\end{figure}
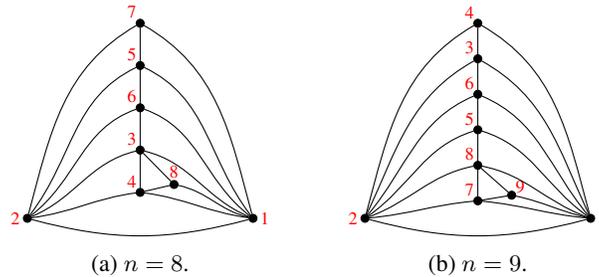

\begin{figure}[!htbp]
 \centering
 \subfloat[Octahedral graph.]{
 \begin{tikzpicture}[%FIGURE 7 (A)
  no/.style = {draw,circle,inner sep=1pt,minimum size=4pt,fill=black},
  every label/.append style={label distance=-2pt,text=red, font=\footnotesize},
  transform shape,
  scale = 0.75
 ]
  \node[no,label=180:5] (a) at (0.00,0.00) {};
  \node[no,label=000:4] (b) at (4.00,0.00) {};
  \node[no,label=090:6] (c) at (2.00,3.46) {};
  \node[no,label=000:2] (d) at (2.66,1.54) {};
  \node[no,label=270:3] (e) at (2.00,0.38) {};
  \node[no,label=180:1] (f) at (1.33,1.54) {};
  \draw (a) to [out=-15,in=195] (b)
        (b) to [out=105,in=-30] (c)
        (c) to [out=210,in=75]  (a);
  \draw (a) -- (f) -- (c) -- (d) -- (b) -- (e) -- (d) -- (f) -- (e) -- (a);
 \end{tikzpicture}
  \label{sfig:octahedral-i}
 }
%  \quad
 \subfloat[Heptahedral graph 34.]{
 \begin{tikzpicture}[%FIGURE 7 (B)
  no/.style = {draw,circle,inner sep=1pt,minimum size=4pt,fill=black},
  every label/.append style={label distance=-2pt,text=red, font=\footnotesize},
  transform shape,
  scale = 0.75
 ]
  \node[no,label=180:3] (a) at (0.00,0.00) {};
  \node[no,label=000:4] (b) at (4.00,0.00) {};
  \node[no,label=090:7] (c) at (2.00,3.46) {};
  \node[no,label=000:6] (d) at (2.66,1.54) {};
  \node[no,label=270:2] (e) at (2.00,0.38) {};
  \node[no,label=180:5] (f) at (1.33,1.54) {};
  \node[no,label=270:1] (g) at (2.00,2.18) {};
  \draw (a) to [out=-15,in=195] (b)
        (b) to [out=105,in=-30] (c)
        (c) to [out=210,in=75]  (a);
  \draw (c) -- (f) -- (a) -- (e) -- (b) -- (d) --
        (f) -- (e) -- (d) -- (c) -- (g) -- (f)
        (d) -- (g);
 \end{tikzpicture}
  \label{sfig:graph34-i}
 }
%  \quad
 \subfloat[Johnson solid skeleton 13.]{
 \begin{tikzpicture}[%FIGURE 7 (C)
  no/.style = {draw,circle,inner sep=1pt,minimum size=4pt,fill=black},
  every label/.append style={label distance=-2pt,text=red, font=\footnotesize},
  transform shape,
  scale = 0.75
 ]
  \node[no,label=180:1] (a) at (0.00,0.00) {};
  \node[no,label=000:3] (b) at (4.00,0.00) {};
  \node[no,label=090:7] (c) at (2.00,3.46) {};
  \node[no,label=000:2] (d) at (2.66,1.54) {};
  \node[no,label=270:6] (e) at (2.00,0.38) {};
  \node[no,label=180:4] (f) at (1.33,1.54) {};
  \node[no,label={[label distance=-4pt]230:5}] (g) at (2.00,1.54) {};
  \draw (a) to [out=-15,in=195] (b)
        (b) to [out=105,in=-30] (c)
        (c) to [out=210,in=75]  (a);
  \draw (c) -- (f) -- (a) -- (e) -- (b) -- (d) --
        (c) -- (g) -- (e) -- (f) -- (g) -- (d) -- (e);
 \end{tikzpicture}
  \label{sfig:johnson13-i}
 }
 \caption{Some graph numberings for which Lemma \ref{lbl-apollonian}-\ref{a-i} does not hold.}
 \label{fig:nonappolonian-i}
\end{figure}
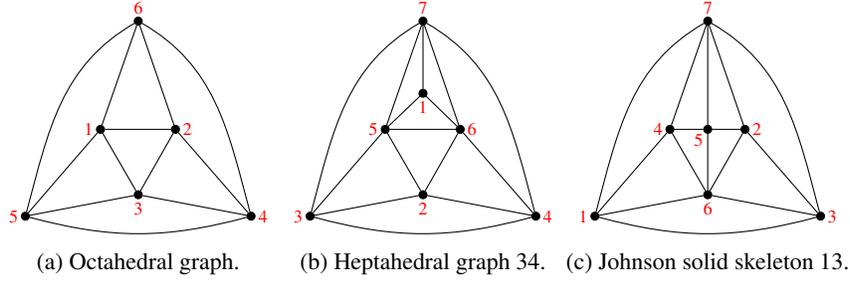

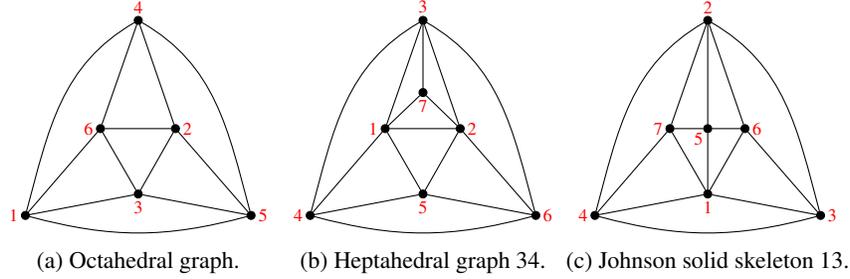
\begin{figure}[!htbp]
 \centering
 \subfloat[Octahedral graph.]{
 \begin{tikzpicture}[%FIGURE 8 (A)
  no/.style = {draw,circle,inner sep=1pt,minimum size=4pt,fill=black},
  every label/.append style={label distance=-2pt,text=red, font=\footnotesize},
  transform shape,
  scale = 0.75
 ]
  \node[no,label=180:1] (a) at (0.00,0.00) {};
  \node[no,label=000:5] (b) at (4.00,0.00) {};
  \node[no,label=090:4] (c) at (2.00,3.46) {};
  \node[no,label=000:2] (d) at (2.66,1.54) {};
  \node[no,label=270:3] (e) at (2.00,0.38) {};
  \node[no,label=180:6] (f) at (1.33,1.54) {};
  \draw (a) to [out=-15,in=195] (b)
        (b) to [out=105,in=-30] (c)
        (c) to [out=210,in=75]  (a);
  \draw (a) -- (f) -- (c) -- (d) -- (b) -- (e) -- (d) -- (f) -- (e) -- (a);
 \end{tikzpicture}
  \label{sfig:octahedral-ii}
 }
%  \quad
 \subfloat[Heptahedral graph 34.]{
 \begin{tikzpicture}[%FIGURE 8 (B)
  no/.style = {draw,circle,inner sep=1pt,minimum size=4pt,fill=black},
  every label/.append style={label distance=-2pt,text=red, font=\footnotesize},
  transform shape,
  scale = 0.75
 ]
  \node[no,label=180:4] (a) at (0.00,0.00) {};
  \node[no,label=000:6] (b) at (4.00,0.00) {};
  \node[no,label=090:3] (c) at (2.00,3.46) {};
  \node[no,label=000:2] (d) at (2.66,1.54) {};
  \node[no,label=270:5] (e) at (2.00,0.38) {};
  \node[no,label=180:1] (f) at (1.33,1.54) {};
  \node[no,label=270:7] (g) at (2.00,2.18) {};
  \draw (a) to [out=-15,in=195] (b)
        (b) to [out=105,in=-30] (c)
        (c) to [out=210,in=75]  (a);
  \draw (c) -- (f) -- (a) -- (e) -- (b) -- (d) --
        (f) -- (e) -- (d) -- (c) -- (g) -- (f)
        (d) -- (g);
 \end{tikzpicture}
  \label{sfig:graph34-ii}
 }
%  \quad
 \subfloat[Johnson solid skeleton 13.]{
 \begin{tikzpicture}[%FIGURE 8 (C)
  no/.style = {draw,circle,inner sep=1pt,minimum size=4pt,fill=black},
  every label/.append style={label distance=-2pt,text=red, font=\footnotesize},
  transform shape,
  scale = 0.75
 ]
  \node[no,label=180:4] (a) at (0.00,0.00) {};
  \node[no,label=000:3] (b) at (4.00,0.00) {};
  \node[no,label=090:2] (c) at (2.00,3.46) {};
  \node[no,label=000:6] (d) at (2.66,1.54) {};
  \node[no,label=270:1] (e) at (2.00,0.38) {};
  \node[no,label=180:7] (f) at (1.33,1.54) {};
  \node[no,label={[label distance=-4pt]230:5}] (g) at (2.00,1.54) {};
  \draw (a) to [out=-15,in=195] (b)
        (b) to [out=105,in=-30] (c)
        (c) to [out=210,in=75]  (a);
  \draw (c) -- (f) -- (a) -- (e) -- (b) -- (d) --
        (c) -- (g) -- (e) -- (f) -- (g) -- (d) -- (e);
 \end{tikzpicture}
  \label{sfig:johnson13-ii}
 }
 \caption{Some graph numberings for which Lemma \ref{lbl-apollonian}-\ref{a-ii} does not hold.}
 \label{fig:nonappolonian-ii}
\end{figure}

The following two lemmas express the minimal and maximal quantities of valid paths for two-dimensional (planar) square grid graphs that have an odd number of rows and columns. It can be seen that the maximum quantity is always exactly three times the minimum quantity.

% _____________ O D D   G R I D S
\begin{lemma}\label{lbl-grid-odd}
If $G_{n\times n}$ is a two-dimensional square grid graph having $n$ rows and $n$ columns, with $n$ odd, then ($i$) $\phi_{\min}(G_{n\times n}) = (n - 1)^2$ and ($ii$) $\phi_{\max}(G_{n\times n}) = 3(n - 1)^2$.
\end{lemma}
\begin{proof}
Let $G_{n\times n}$ be a two-dimensional (planar) square grid graph having $n$ rows and $n$ columns, where $n = 2k + 1$ for some integer $k$. 

\begin{enumerate}[leftmargin=*,label=(\emph{\roman*})]
\item The following proof rests on constructing a minimal numbering of $G_{n\times n}$ and showing that it induces $(n-1)^2$ valid paths. Let the vertex of $G_{n\times n}$ occupying the $i$th row and $j$th column in a planar embedding of $G_{n\times n}$ be denoted as the ordered pair $(i, j)$, for $i, j = 1,\dots,n$. Define a numbering $\pi_{OddMin}$, that numbers the vertices of $G_{n\times n}$ row-wise from the top left vertex to the bottom-right vertex, i.e., in the order $(1,1), (1,2), \dots, (1,n),\allowbreak  (2,1), (2,2), \dots, (2,n),\dots, (n,1), (n,2), \dots, (n,n)$; as $1, 2, \dots, n^2$, respectively. Let $S(G_{n\times n})$ denote the vertex set containing all the vertices of $G_{n\times n}$ except those of the right-most column and the bottom row, i.e., $\{(1,1), (1,2), \dots,(1,n-1), (2,1),\allowbreak (2,2), \dots, (2,n-1), \dots, (n-1,1), (n-1,2), \dots, (n-1,n-1)\}$. Note that the set $S(G_{n\times n})$ contains $4k^2$ vertices. (There are $(2k-1)^2$ internal vertices of $G_{n\times n}$ (all of degree $4$), $4k-2$ vertices on the perimeter of $G_{n\times n}$ of degree $3$, and exactly one vertex $(1,1)$ of degree $2$.)

By the nature of $\pi_{OddMin}$, each vertex in $S(G_{n\times n})$ is the middle vertex of exactly one valid path in $G_{n\times n}$. Hence, $\pi_{OddMin}$ induces $4k^2 = (n-1)^2$ valid paths based on the vertices of $S(G_{n\times n})$. However, $\pi_{OddMin}$ does not induce any further valid paths (from the vertices excluded from 
$S(G)$). As $\pi$ is a minimal numbering, the result follows.  
\item The proof is based on constructing a maximal numbering of $G_{n\times n}$ and showing that it induces $3(n - 1)^2$ valid paths. Let $S_O(G_{n\times n})$ be the set of vertices of $G_{n\times n}$ occupying the $i$-th row and $j$-th column, for $i + j$ odd. Note that $S_O(G_{n\times n})$ contains $2(k+1)k$ vertices. (There are $2k(k-1)$ internal vertices of $G$ (all of degree $4$) and $4k$ vertices on the perimeter of $G$ of degree $3$.) Define a numbering $\pi_{OddMax}$, that numbers the elements of $S_O(G_{n\times n})$ arbitrarily as $1,2,\dots,2k(k+1)$. Further, $\pi_{OddMax}$ arbitrarily numbers the remaining vertices of $G_{n\times n}$ (with $i + j$ even) with the rest of the numbers, being $2k^2+2k+1,\dots,n^2$. Let $u$ be any vertex of $S_O(G_{n\times n})$. By the nature of $\pi_{OddMax}$, all the 2-paths having $u$ as a middle vertex are valid. If $u$ is of degree $4$ ($3$) it is the middle vertex of $6$ ($3$) 2-paths. Hence, $\pi_{OddMax}$ creates a total of $6(2k(k-1)) + 3(4k) = 3(n-1)^2$ valid paths based on the vertices of $S_O(G_{n\times n})$. $\pi_{OddMax}$ does not induce any further valid paths (from the vertices excluded from $S_O(G_{n\times n})$). As $\pi_{OddMax}$ is a maximal numbering, the result follows.
\end{enumerate}
\end{proof}

\begin{table}[!htbp]
$$
 \begin{array}{ccccccccccc}
  \hline
 k  &\ &   n  &\ &  2k(k-1)\cdot 6  &\ &  4k\cdot 3 &\ &  \phi_{\min}(G_{n\times n})   &\ &  \phi_{\max}(G_{n\times n})   \\
\hline
1 && 3 &&   0\cdot 6  &&   4\cdot 3 && 4 && 12 \\
2 && 5 &&   4\cdot 6  &&   8\cdot 3 && 16 && 48 \\
3 && 7 &&  12\cdot 6  &&  12\cdot 3 && 36 && 108 \\
4 && 9 &&  24\cdot 6  &&  16\cdot 3 && 64 && 192 \\
5 && 11 &&  40\cdot 6  &&  20\cdot 3 && 100 && 300 \\
\vdots\\
 k  &&   2k+1    &&  2k(k-1)\cdot 6  &&  4k\cdot 3 &&  4k^2=(n-1)^2  &&  12k^2 =3(n-1)^2 \\
 %   &&           &&                  &&            &&  (n-1)^2       &&  3(n-1)^2 \\
\hline
 \end{array}
$$
\caption{Optimal numberings of square grid graphs with an odd number of rows and columns.}
\label{tab:grids-odd}
\end{table}

Optimal numberings for grid graphs of odd degree are shown in Table \ref{tab:grids-odd}. Recall that each graph in Table \ref{tab:grids-odd} is of order $n$. Thus, when graph $G$ is a grid graph it is assumed that $n$ is a perfect square and the number of rows and columns are both $\sqrt{n}$.  This is at variance to Lemma \ref{lbl-grid-odd} and Table \ref{tab:grids-odd} where the grid graphs are of order $n^2$. 

A summary of the above findings are given in Table \ref{tab:complexities}. As can be seen from the table, for the totality of graph classes studied, apparently there is no recognisable general pattern for a numbering that leads to the optimisation (minimisation  or maximisation) of the quantity of valid paths. Furthermore, each graph class in the table has its own individual pattern. These observations strengthen the conjecture that such a pattern effectively does not exist.
\begin{table}[!htbp]
$$
\begin{array}{lcc}
\hline\\[-10pt]
G & \phi_{\min}(G) & \phi_{\max}(G)\\[2pt]
\hline\\[-6pt]
P_n\>(\text{Path}) & 0 & \lceil\frac{n}{2}\rceil-1\\[6pt]
T_n\>(\text{Tree}) & 0 & \binom{n-1}{2}\\[6pt]
C_n\>(\text{Cycle}) & 1 & \left\lceil\frac{n-1}{2}\right\rceil\\[6pt]
W_n \>(\text{Wheel}) & n & \binom{n-1}{2}+\left\lceil\frac{n-1}{2}\right\rceil \\[6pt]
K_{p,q}\>(\text{Complete bipartite)} &  (3q-p-1)(p^2-p)/6 & p\binom{q}{2} \\[6pt]
K_n \>(\text{Complete})& \sum^{n-1}_{i = 2}\binom{i}{2} & \sum^{n-1}_{i = 2}\binom{i}{2}\\[6pt]
A_n\>(\text{Apollonian, $n=  5$}) & 3n-8 & 3n-5 \\[6pt]
A_n\>(\text{Apollonian, $n=  6$}) & 3n-8 & 3n-1 \\[6pt]
A_n\>(\text{Apollonian, $n \geqslant 7$}) & 3n-8 & [5n-12, 5n-11] \\[6pt]
G_{n,n}\>(\text{Square grid, $n$ odd}) & (n-1)^2 & 3(n-1)^2\\[6pt]
\hline
\end{array}
$$
\caption{$\phi_{\min}(G)$ and $\phi_{\max}(G)$  for the optimal numberings of some classes of graphs.}
\label{tab:complexities}
\end{table}

%--------------------------------------------------- S E C T I O N -
\section{Conclusion}
\label{sec:conclusion}

The current work introduces a new kind of graph numbering, which opens a new field of study in graph labelling. Two related optimising problems were defined and some findings concerning particular classes of graphs were presented. The results showed that apparently there does not exist a particular numbering that is applicable to more than one of the classes of graphs under analysis, which requires the development of a specific numbering for each graph class studied.

Therefore, the proposed problems are nontrivial, leading to Conjectures \ref{conjmin} and \ref{conjmax}, given next. They state that both MIN-VP and MAX-VP belong to the $\mathcal{NP}$-complete class of problems:

\begin{conjecture}\label{conjmin} Deciding, for a given graph $G$ and a given integer $k$, whether or not $G$ has a numbering that induces at most $k$ valid paths is an $\mathcal{NP}$-complete decision problem. 
\end{conjecture}

\begin{conjecture}\label{conjmax} Deciding, for a given graph $G$ and a given integer $k$, whether or not $G$ has a numbering that induces at least $k$ valid paths is an $\mathcal{NP}$-complete decision problem. 
\end{conjecture}

%--------------------------------------------- R E F E R E N C E S -
\bibliographystyle{abbrvnat}
\bibliography{graph-labeling}
\end{document}